\documentclass[10pt,a4paper]{article}
\usepackage[utf8]{inputenc}
\usepackage{amsmath}
\usepackage{amsfonts}
\usepackage{amssymb}
\usepackage{amsmath,amsxtra,amssymb,latexsym,amscd,amsthm,amsfonts}
\newtheorem{theorem}{Theorem}
\newtheorem{cor}[theorem]{Corollary}

\newtheorem{lem}[theorem]{Lemma}

\newtheorem{rem}[theorem]{Remark}

\DeclareMathOperator{\divv}{div}

\begin{document}

\title{\bf Halfspace type theorems for self-shrinkers in arbitrary codimension}
\author{\bf  Doan The Hieu and Nguyen Thi My Duyen\\
 Department of Mathematics\\
 College of Education, Hue University, Hue, Vietnam\\
 \\
 dthieu@hueuni.edu.vn, ntmduyen@hueuni.edu.vn}

\maketitle
\begin{abstract}
  In this paper, we generalize some  halfspace type theorems for  self-shrinkers of codimension 1 to the case of arbitrary codimension. 
 \end{abstract}

\vskip .5cm
\section{Introduction}

The halfspace theorem says that  {\sl ``There is no non-planar, complete, minimal surface properly immersed in a halfspace of $\mathbb R^3.$''}
The theorem is due to Hoffman and Meeks. In fact they proved a stronger version, the strong halfspace theorem,  {\sl ``Two disjoint complete properly immersed minimal surfaces in $\mathbb R^3$ are planes''} (see \cite{home}).

The halfspace theorem  is essentially a three-dimensional one. In $\mathbb R^n, n> 3,$ the halfspace theorem is false because there are minimal Catenoids with bounded height.

Many generalizations of the theorem have been made by several authors, see \cite{daha}, \cite{damero}, \cite{ma}, \cite{nesa}, \cite{roro}, \cite{roscsp} and references therein. 

The first halfspace theorem for self-shrinker in codimension 1  was proved in \cite{piri} based on the weighted parabolicity of self-shrinkers. A similar result in a more general setting was proved recently in \cite{huparo}.
\begin{theorem}[Theorem 3 in \cite{piri}; Theorem 1.1 in \cite{caes}]\label{theo1.1} 
Let $P$ be a hyperplane passing through the origin. The only properly
immersed self-shrinker contained in one of the closed halfspace determined by $P$ is $\Sigma = P.$
\end{theorem}
In contrast with  the case of minimal surfaces, the halfspace theorem for self-shrinkers holds true in any dimension. Moreover, one can consider a type of  halfspace theorems for self-shrinker containing inside or outside a hypercylinder. 

In 2016, Cavalcante and  Espinar \cite{caes} showed some halfspace type theorems for  self-shrinkers of codimension 1 including Theorem \ref{theo1.1} with a different proof. 
\begin{theorem} [Theorem 1.2 in \cite{caes}]
The only complete self-shrinker properly immersed in a closed cylinder
$\overline{B^{k+1}(R)} \times \mathbb R^{n-k}\subset \mathbb R^{n+1},$ for some $k\in\{1, \dots , n\}$ and radius $R, R\le\sqrt{2k},$ is the cylinder
$S^k(\sqrt{2k})\times \mathbb R^{n-k}.$
\end{theorem}
\begin{theorem} [Theorem 1.3 in \cite{caes}]
The only complete self-shrinker properly immersed in an exterior closed
cylinder $\overline{E^{k+1}(R)}\times\mathbb R^{n-k} \subset\mathbb R^{n+1},$ for some $k \in \{1, \ldots , n\}$ and radius $R, R\ge \sqrt{2k},$ is the
cylinder $S^k(\sqrt{2k}) \times\mathbb R^{n-k}.$ Here $E^{k+1}(R)=R^{k+1}-\overline{B^{k+1}(R)}.$
\end{theorem}
In  2018, Vieira and  Zhou \cite{vizh} proved similar results, where spheres or balls center at the origin are replaced by ones with arbitrary centers and suitable radius. Recently, Imper, Pigola and Rimoldi \cite{impiri} recovered Cavalcante and Espinar's results with  short proofs by using potential theoretic arguments.

 The paper aims  to generalize the above haflspace  type results for codimension 1 self-shrinkers  to the case of arbitrary codimension. The first step in our approach is somewhat similar to the one in \cite{impiri} for codimension 1 but the use of maximal principle for weighted superharmonic functions together with the weighted parabolicity of self-shrinkers is replaced by an application of a divergence theorem (Theorem \ref{theo6}).  In fact our proofs recovered some key formulas 
 (\ref{eq5}), (\ref{eq7}), (\ref{eq9}) that are due Colding-Minicozzi  \cite{comi3}  for the case of codimension 1.   Arezzo-Sun  \cite{aesu} observed that these formulas are also  true for the case of arbitrary codimension.
 
We would like to thank  Vieire, Rimoldi, Rosales for introducing us to their interesting  works and the others for helpful comments and suggestions.
 \section{Preliminaries}
 In this paper, we use the following notations
\begin{enumerate}
\item $B^k(a, R),$  the $k$-ball with center $a$ and radius $R;$
 \item $E^k(a, R)=\mathbb R^k-\overline{B^k(a, R)},$ the complement of  $\overline{B^k(a,R)};$
 \item $S^k(a, R),$ the $k$-sphere with center $a$ and radius $R;$
 \item  $\overline {A},$ the closure of the set $A.$
 \end{enumerate}
  For simple, when the center of spheres or balls is the origin we write $B^k(R), E^k(R), S^k(R).$ 
 \subsection{Self-shrinkers}

An $n$-dimensional submanifold $\Sigma$ immersed  in $\mathbb R^m, m>n, $ is called  a self-shrinker for
the mean curvature flow (MCF), if 
\begin{equation}\label{se}
{\bf H} = -\frac 12X^N,
\end{equation}  
where ${\bf H}$ is the mean curvature vector of $\Sigma, X$  is the position vector, and $X^N$ denotes  the
normal part of $X.$
 
 Self-shrinkers are self-similar solutions to MCF and play an important role in the study of singularities of the flow. For more information about self-shrinkers  as well as singularities, we refer
the readers to \cite{comi1}, \cite{comi2}, \cite{hu}, \cite{hu1}.
 
A complete self-shrinker $\Sigma^{n}$ in $\mathbb{R}^{m}$ is said to have  polynomial
volume growth if there exist constants $C_1$ and $d_1$ such that for all $R\ge 1,$ there holds
\begin{equation}\label{poly}
{\rm Vol} (B^m(R)\cap \Sigma)\le C_1R^{d_1}.
\end{equation}
In 2013, Cheng-Zhou  \cite{chzh} and  Ding-Xin \cite {dixi}, proved that 

``{\sl
A complete non-compact properly immersed
self-shrinker $\Sigma^{n}$ in $\mathbb{R}^{m}, m>n,$ has Euclidean volume growth at most, i.e.
$${\rm Vol}(B^m(R)\cap \Sigma)\le CR^{n}$$ 
for $R\geq 1.$}''

 \subsection{Some typical examples}
 It is not hard to verify all of the followings are $n$-dimensional complete self-shrinkers in $\mathbb R^m.$
 \begin{enumerate}
 \item An $n$-plane passing through the origin.
 \item $S^n(\sqrt{2n})\subset \mathbb R^{n+1}.$ 
 \item The cylinder $S^k(\sqrt{2k})\times\mathbb R^{n-k}\subset \mathbb R^{n+1},  0<k<n. $
 \item $S^{n_1}(\sqrt{2n_1})\times S^{n_2}(\sqrt{2n_2})\times\ldots \times S^{n_k}(\sqrt{2n_k})\subset \mathbb R^{n+1},$  $ n_1+n_2+\ldots n_k=n.$
  \item $S^{n_1}(\sqrt{2n_1})\times S^{n_2}(\sqrt{2n_2})\times\ldots \times S^{n_k}(\sqrt{2n_k})\times\mathbb R^p\subset \mathbb R^{n+1},$  $ p\ge1$ and $n_1+n_2+\ldots n_k+p=n.$
 \item $n$-dimensional complete minimal submanifolds of the sphere $S^{m-1}(\sqrt{2n})$ (see Theorem 4.1 in \cite {aesu} or subsection 1.4 in \cite{xi}).
 \end{enumerate} 
 
 For some more well-known results about complete self-shrinkers, we refer the readers to \cite{comi3}, \cite{hu1}, \cite{piri} for the case of codimension 1 and \cite{chpe}, \cite{sm} for the case of arbitrary codimension.

\subsection{Some  calculations} 
  In this subsection, we calculate  the surface divergence of some vector fields that will be used in the proofs of the main results. The calculations are straightforward, but for the sake of completeness we present them here.
 
Let $e_1,e_2, \ldots, e_m$ be the coordinate vector fields for $\mathbb R^m,$ $ \Sigma^n$ be a complete self-shrinker in $\mathbb R^m,  \{E_1,E_2,\ldots, E_n\}$ be an orthonormal basis for $T_X\Sigma, X=\sum_{i=1}^{m}x_ie_i$ be the position vector field and  $u=\sum_{i=1}^{k+1}x_ie_i, k\le m-1.$
 We have the following lemma.
\begin{lem}  \label{lem1}  
\begin{enumerate}
\item 
\begin{equation}\label{xt}
\divv_\Sigma X^T=n-\frac12 |X^N|^2.
\end{equation}
\item
\begin{equation} \label{et}
\divv_\Sigma e_l^T=-\frac 12\langle X, e_l^N\rangle,\ l=1,2,\ldots, m.
\end{equation}
\item
\begin{equation}\label{xet}
\divv_\Sigma x_le_l^T=|e_l^T|^2-\frac12x_l\langle X, e_l^N\rangle,\ l=1,2,\ldots, m.
\end{equation}
\item 
\begin{equation}\label{ut}
\divv_\Sigma u^T=(k+1) -\frac12|u^N|^2-\sum_{i=1}^{k+1}|e_i^N|^2.
\end{equation}
\item
\begin{equation}\label{uu}
\divv_\Sigma \frac 1{|u|}u^T=\frac 1{|u|}\left[k -\frac12|u^N|^2-\sum_{i=1}^{k+1}|e_i^N|^2+\frac {|u^N|^2}{|u|^2}\right].
\end{equation}

\end{enumerate}
\end{lem}
\begin{proof}
 We use the summation convention.
\begin{enumerate}
\item We have
$$\divv_\Sigma X=n,$$
and
\begin{align*}
\divv_\Sigma X^N&=\langle E_i, \nabla_{E_i}X^N\rangle=\nabla_{E_i}\langle E_i, X^N\rangle-\langle\nabla_{E_i}E_i, X^N\rangle\\
&=\nabla_{E_i}(0)-\langle (\nabla_{E_i}E_i)^N, X\rangle=-\langle {\bf H}, X\rangle=\frac12 |X^N|^2.
\end{align*}
Therefore, 
\begin{align*}
\divv_\Sigma X^T=n-\frac12 |X^N|^2.
\end{align*}
\item 
\begin{align*}
\divv_\Sigma e_l^T&=\divv_\Sigma e_l-\divv_\Sigma e_l^N= 0- \langle E_i,\nabla_{E_i}e_l^N\rangle\\
&=\langle \nabla_{E_i}E_i, e_l^N\rangle=\langle (\nabla_{E_i}E_i)^N, e_l\rangle=\langle {\bf H}, e_l\rangle\\
&=-\frac12\langle X, e_l^N\rangle.
\end{align*}
\item
\begin{align*}
\divv_\Sigma x_le_l^T&=\divv_\Sigma x_le_l-\divv_\Sigma x_le_l^N= |e_l^T|^2- \langle E_i,\nabla_{E_i}x_le_k^N\rangle\\
&=|e_l^T|^2+\langle(\nabla_{E_i}E_i)^N, x_le_l\rangle=|e_l^T|^2+\langle{\bf H}, x_le_l\rangle\\
&=|e_l^T|^2-\frac12x_l\langle X, e_l^N\rangle.
\end{align*}

\item For $v\in T_p\Sigma,$ 
$$\nabla_vu=\pi_1(v)=\langle v, e_1\rangle e_1+\langle v, e_2\rangle e_2+\ldots, \langle v, e_{k+1}\rangle e_{k+1}.$$
We have
\begin{align*}
\divv_\Sigma(u)&=\langle E_i,\nabla_{E_i}u\rangle=\sum_{j=1}^{k+1}\sum_{i=1}^n\langle E_i, e_j\rangle^2\\
&=\sum_{j=1}^{k+1}|e_j^T|^2=(k+1)-\sum_{j=1}^{k+1}|e_j^N|^2,\\
\end{align*}
and
\begin{align*}
\divv_\Sigma u^N&=\langle E_i, \nabla_{E_i}u^N\rangle=\nabla_{E_i}\langle E_i, u^N\rangle-\langle \nabla_{E_i}E_i, u^N\rangle\\
&=\nabla_{E_i}(0)-\langle(\nabla_{E_i}E_i)^N, u\rangle=-\langle {\bf H}, u\rangle=\frac12|u^N|^2.
\end{align*}
Therefore,
\begin{align*}
\divv_\Sigma u^T=(k+1) -\frac12|u^N|^2-\sum_{i=1}^{k+1}|e_i^N|^2.
\end{align*}
\item
\begin{align*}
\divv_\Sigma \frac 1{|u|}u^T&=\langle \nabla_\Sigma\frac 1{|u|}, u^T\rangle+\frac 1{|u|}\divv u^T\\
&=-\frac {|u^T|^2}{|u|^3}+\frac 1{|u|}[(k+1) -\frac12|u^N|^2-\sum_{i=1}^{k+1}|e_i^N|^2]\\
&=\frac 1{|u|}\left[k -\frac12|u^N|^2-\sum_{i=1}^{k+1}|e_i^N|^2+\frac {|u^N|^2}{|u|^2}\right].
\end{align*}
\end{enumerate}
\end{proof}

\section{Results} 
In this section, $\Sigma$ is assumed to be an $n$-dimensional complete (without boundary)  self-shrinker properly immersed  in $\mathbb R^m, m>n.$ 

The condition of polynomial volume growth is essential for using an integral formula that is similar to the generalized divergence theorem for compact manifolds. We have the following theorem.

 \begin{theorem}\label{theo6}
Let $F$ be a smooth tangent vector field on $\Sigma.$ For every $X\in\Sigma,$ if $ |\divv_\Sigma F(X)|\le C_2|X|^{d_2},$ where $C_2$ is a positive constant and $d_2$ is a positive integer, then 
\begin{equation}
\int_\Sigma \divv_\Sigma (e^{-\frac{X^2}4}F)dV=0.
\end{equation}
\end{theorem}

\begin{proof}

We only need to prove for the case $\Sigma$ is non-compact. Since $\Sigma$ is proper, $\partial(B_R\cap\Sigma)\ne\emptyset$ when $R$ is large enough.
Since $F$ is tangent to $\Sigma,$ the generalized divergence theorem  for $e^{-\frac{X^2}4}F$ yields
$$\int_{B_R\cap\Sigma} \divv_\Sigma (e^{-\frac{X^2}4}F)dV=e^{-\frac{R^2}4}\int_{\partial(B_R\cap\Sigma)}\left\langle F,\nu\right\rangle dA.$$
Taking the limit when $R\rightarrow\infty,$ the theorem is proved because
\begin{align*}
\lim_{R\rightarrow\infty}e^{-\frac{R^2}4}\left |\int_{\partial(B_R\cap\Sigma)}\left\langle F,\nu\right\rangle dA\right |&=\lim_{R\rightarrow\infty}e^{-\frac{R^2}4}\left|\int_{B_R\cap\Sigma} \divv_\Sigma F dV\right|\\
&\le \lim_{R\rightarrow\infty}e^{-\frac{R^2}4}C_2|X|^{d_2}\int_{B_R\cap\Sigma} dV\\
&\le \lim_{R\rightarrow\infty}e^{-\frac{R^2}4}C_1C_2R^{d_1+d_2}=0.
\end{align*}

\end{proof}
Applying Theorem \ref{theo6} with suitable choices of tangent vector fields $F,$ we obtain the main results of the paper.
 \subsection{Half space type result w.r.t. hyperplanes}
The following theorem says that  $\Sigma$ intersects every hyperplane  passing through the origin.
\begin{theorem} 
\label{theoP}
Let $P$ be a hyperplane passing through the origin. If $\Sigma$ lies  in a closed halfspace determined by $P,$ then  $\Sigma \subset P.$ 
\end{theorem}

\begin{proof}Without loss of generality, we can suppose that $P$ is the hyperplane $x_m=0$ and $\Sigma$ is in the closed half space $\{(x_1, x_2,\ldots, x_m):x_m\ge0\}.$

By (\ref{et}),
\begin{align*}
\divv_\Sigma (e^{-\frac{X^2}4}e_m^T)&=e^{-\frac{X^2}4}\divv_\Sigma e_m^T- e^{-\frac{X^2}4} \frac 12\langle X, e_m^T\rangle \\ \nonumber
&=-\frac 12 e^{-\frac{X^2}4}\left[\langle X, e_m^N\rangle  +\langle X, e_m^T\rangle        \right]\\
&=-\frac12 e^{-\frac{X^2}4}x_m.
\end{align*}
Then Theorem \ref{theo6} applying for $F=e_m^T$ yields (see \cite{comi3} for the case of codimension 1, also see \cite{aesu})
\begin{equation} \label{eq5}
\int_\Sigma e^{-\frac{X^2}4}x_mdV=0.
\end{equation}
Therefore, $x_m=0,$ i.e. $\Sigma\subset P.$

\end{proof}
\begin{rem} 
 If $n=m-1,$ then $\Sigma= P$ (\cite{piri}, Theorem 3 ;  \cite{caes}, Theorem 1.1).
\end{rem}
\begin{cor} \label{co10}
If there exist $m-n$ orthonormal vectors $v_1, v_2,\ldots, v_{m-n}$ such that for $i=1,2,\ldots, m-n, \langle X, v_i\rangle$  does not change sign,  then $\Sigma$ is an $n$-plane passing through the origin.
\end{cor}
\begin{proof}
Without loss of generality, we can assume that $v_i=e_{n+i}$ if $\langle X, v_i\rangle\ge 0$ and  $v_i=-e_{n+i}$ if $\langle X, v_i\rangle\le 0.$ The assumption guarantees that $\Sigma$ is in the closed halfspace $\{(x_1,x_2,\ldots, x_m): x_{n+i}\ge 0, i=1,2,\ldots, m-n\}.$  The proof is then followed by applying Theorem \ref{theoP} in turn for $v_1, v_2,\ldots, v_{m-n}.$
\end{proof}
Based on the  Bernstein result  for self-shrinkers of codimension 1, 
 {\sl ``An entire graphic self-shrinker must be a hyperplane passing through the origin''}  (see \cite {echu}, \cite{wa}, \cite{hi}), and
with the same argument as in the proof of Corollary \ref{co10},  we have the following.
\begin{cor} [A Bernstein type theorem] 
Let  $F:\mathbb R^n\rightarrow \mathbb R^{m-n}, F({\bf x})=(f_1({\bf x}), f_2({\bf x}),\ldots, f_{m-n}({\bf x}))$ be a smooth function and $\Sigma=\{({\bf x}, F({\bf x})): {\bf x}\in \mathbb R^n\}$ be its graph. If there exist at least $(m-n-1)$ functions $f_i$ that do not change sign, then $\Sigma$ is an $n$-plane passing through the origin. 
\end{cor}

\subsection{Self-shrinkers inside or outside a ball}

The following theorem says that a complete properly immersed  self-shrinker $\Sigma^n$ and $S^{m-1}(\sqrt{2n})$ must be intersected.
\begin{theorem}\label{theoball}
If $\Sigma\subset \overline{E^m(\sqrt{2n})}$ or  $\Sigma\subset \overline {B^{m}(\sqrt{2n})},$ then $\Sigma$ is compact and $\Sigma\subset S^{m-1}(\sqrt{2n}),$ i.e. $\Sigma$ is a minimal submanifold of $S^{m-1}(\sqrt{2n}).$  Moreover, if $n=m-1,$ then $\Sigma= S^{n}(\sqrt{2n}).$ 
\end{theorem}
\begin{proof}
By (\ref{xt}),
\begin{align*}
 \divv_\Sigma (e^{-\frac{X^2}4}X^T)&=e^{-\frac{X^2}4}\divv_\Sigma X^T-e^{-\frac{X^2}4}\langle \frac 12X, X^T\rangle \\ \nonumber
&= e^{-\frac{X^2}4}(n-\frac 12|X^N|^2)- e^{-\frac{X^2}4}\frac 12 |X^T|^2\\
&=e^{-\frac{X^2}4}(n-\frac 12|X|^2).
\end{align*}

Applying Theorem \ref{theo6} with $F=X^T$ (see \cite{comi3} for the case of codimension 1, also see \cite{aesu}),
\begin{equation}\label{eq7}
\int_\Sigma e^{-\frac{X^2}4}(n-\frac 12|X|^2)dV=0.
\end{equation}

If $\Sigma\subset \overline{E^m(\sqrt{2n})}$  ($\Sigma\subset \overline {B^{m}(\sqrt{2n})}$), then $2n-|X|^2\le 0$ ($2n-|X|^2\ge 0$). By (\ref{eq7}), it follows that  $2n-|X|^2= 0,$ i.e. $\Sigma\subset S^{m-1}(\sqrt{2n}).$ Since $\Sigma$ is proper, it must be compact.

The case of $n=m-1$ is obvious.
\end{proof}
The following theorem can be seen as an arbitrary codimension version of Theorem 1 in \cite{vizh}. Here the proof is  also applied for the case of self-shrinkers are outside of spheres.  
\begin{theorem}\label{theovizh1}
\begin{enumerate}
\item Any complete self-shrinker $\Sigma^n$ properly immersed in $\mathbb R^m, m>n,$ intersects all
members of the collection $ C$ given by
$$C := \{S^{m-1}(a, \sqrt{2n + |a|^2}): a\ {\text{ is a  vector in}}\ \mathbb R^m\}.$$
 \item If the $\Sigma$ lies in $\overline {B^m(a,\sqrt{2n + |a|^2})}$ or in $\mathbb R^m-B^m(a,\sqrt{2n + |a|^2})$ then $\Sigma\subset  S^{m-1}(a,\sqrt{2n + |a|^2}).$  Moreover,  if $n=m-1,$ then $\Sigma$ is the sphere $S^{n}(\sqrt{2n}).$
\end{enumerate}
\end{theorem}
\begin{proof}
From (\ref{eq5}), it follows that
\begin{equation} \label{eq51}
\int_\Sigma e^{-\frac{X^2}4}\langle X, a\rangle dV=0.
\end{equation}
Therefore, (\ref{eq7}) and (\ref{eq51}) yields
\begin{equation} \label{eq71}
\int_\Sigma e^{-\frac{X^2}4}(|X-a|^2-(2n+|a|^2))dV=0.
\end{equation}
The theorem is proved easily by some arguments as in the proof of Theorem \ref{theoball}. Note that, for codimension 1 case, the sphere  $S^{n}(a,\sqrt{2n + |a|^2})$ is a self-shrinker if and only if $a=0.$
\end{proof}
\begin{rem}
Theorem 5.1  in \cite{gipa} shows another version of Theorem \ref{theoball}, where self-shrinkers are assumed to be parabolic instead of proper. And a different proof of Theorem \ref{theoball}, stated in terms of $\lambda$-self-shrinkers,  was also done in \cite{gipa} (Theorem 6.3).
\end{rem}

  \subsection{Half space type results w. r. t.  cylinders}
  
  \begin{theorem}  [Self-shrinker inside a hypercylinder]
Let $k \in \{m-n, . . . , m-2\}, p=m-k-1$ and $R=\sqrt{2(n-p)}.$ If  $\Sigma$ is inside the closed cylinder $\overline{B^{k+1}(R)}\times \mathbb R^p,$ then $\Sigma\subset S^{k}(R)\times \mathbb R^p .$ 
\end{theorem} 
\begin{proof}
By (\ref{xet})
\begin{align*}
\divv_\Sigma (e^{-\frac{X^2}4}x_ie_i^T)&=e^{-\frac{X^2}4}\left[\divv_\Sigma (x_ie_i^T)-\frac 12\langle X, x_ie_i^T\rangle\right]\\
&=e^{-\frac{X^2}4}\left[|e_i^T|^2-\frac 12x_i\langle X, e_i^N\rangle  -\frac 12x_i\langle X, e_i^T\rangle\right]\\
&=e^{-\frac{X^2}4}\left[|e_i^T|^2-\frac 12x_i^2\right].
\end{align*}

Applying Theorem \ref{theo6} with $F=x_ie_i^T, $  we have (see \cite{comi3} for the case of codimension 1, also see \cite{aesu})
\begin{equation}\label{eq9}
\int_\Sigma e^{-\frac{X^2}4}x_i^2dV=2\int_\Sigma e^{-\frac{X^2}4}|e_i^T|^2dV.
\end{equation}

Let $\{e_1,e_2,\ldots, e_m\}$ be the standard basis in $\mathbb R^m,$ where $\{e_1, e_2, \ldots, e_{k+1}\}\subset\mathbb R^{k+1}$ and $\{e_{k+2}, e_{k+3}, \ldots, e_{m}\}\subset\mathbb R^{p}.$
Denote $X=(u, v),$ where $u\in \mathbb R^{k+1}, v\in\mathbb R^p.$

By (\ref{eq7}) and (\ref{eq9}), we get
\begin{align*}
\int_\Sigma e^{-\frac{X^2}4}\left[|X|^2-2n-\sum_{i=k+2}^m x_i^2\right]dV&=\int_\Sigma e^{-\frac{X^2}4}[|u|^2-2n]dV\\
&=-2\int_\Sigma e^{-\frac{X^2}4}\sum_{i=k+2}^m  |e_i^T|^2dV.
\end{align*}
Since $|e_i^T|^2=1-|e_i^N|^2,$ it follows that
$$
\int_\Sigma e^{-\frac{X^2}4}\left[|u|^2-R^2\right]dV=2\int_\Sigma e^{-\frac{X^2}4}\sum_{i=k+2}^m  |e_i^N|^2dV\ge 0.
$$
The assumption that $\Sigma$ is inside the closed cylinder $\overline{B^{k+1}(R)}\times \mathbb R^p,$  means
$$|u|^2-R^2\le 0.$$
Therefore, 
$$|u|^2-R^2=0,$$
 i.e. $\Sigma\subset S^{k}(R)\times \mathbb R^p .$ 
 
\end{proof}
 \begin{rem}
 \begin{enumerate}
 \item We see in the above proof that  $e_i^N=0,$ i.e. $e_i=e_i^T, i=k+2,\ldots, m.$ Therefore, $\Sigma=\Gamma\times \mathbb R^p,$ where $\Gamma\subset S^{k}$ is an $(n-p)$-dimensional self-shrinker, i.e. an $(n-p)$-dimensional  minimal submanifold of $S^{k}.$ 

\item If $n=m-1,$ then $\Sigma=S^{k}(\sqrt{2k})\times \mathbb R^{n-k}$ (\cite{caes}, Theorem 1.2).
  \end{enumerate}
 \end{rem}
\vskip .3cm
\begin{theorem}[Self-shrinker outside a  hypercylinder]
Let $k \in \{1, . . . , n\}.$  If $\Sigma$ is contained in $\overline{E^{k+1}(\sqrt{2k})}\times \mathbb R^{m-k-1},$ then $\Sigma\subset S^{k}(\sqrt{2k})\times \mathbb R^{m-k-1}.$
\end{theorem} 

{\sl Proof.} Let $u=\sum_{i=1}^{k+1}x_ie_i.$ 
By (\ref{uu})
\begin{align*}\label{eq11}
\divv_\Sigma (e^{-\frac{X^2}4}\frac{1}{|u|}u^T) &= \left[e^{-\frac{X^2}4}\divv_\Sigma(\frac{1}{|u|}u^T)-\frac12\langle X, \frac{u^T}{|u|}\rangle\right] \\ 
&= e^{-\frac{X^2}4}\frac 1{|u|}\left[k -\frac12|u|^2-\sum_{i=1}^{k+1}|e_i^N|^2+\frac {|u^N|^2}{|u|^2}\right].
\end{align*}
It is not hard to check that
$$\sum_{i=1}^{k+1}|e_i^N|^2\ge \frac {|u^N|^2}{|u|^2}.$$
Indeed, we have
\begin{align*}
|u^N|^2&=\left|\sum_{i=1}^{k+1}x_ie_i^N\right|^2=\sum_{i=1}^{k+1}x_i^2|e_i^N|^2+2\sum_{i\ne j}x_ix_j\langle e_i^N, e_j^N\rangle\\
&\le \sum_{i=1}^{k+1}x_i^2|e_i^N|^2+\sum_{i\ne j}x_i^2|e_j^N|^2\\
&\le\left (\sum_{i=1}^{k+1}x_i^2\right)\left(\sum_{i=1}^{k+1}|e_i^N|^2\right)= |u|^2\left(\sum_{i=1}^{k+1}|e_i^N|^2\right).  \\
\end{align*}
Applying Theorem \ref{theo6} with $\displaystyle F=\frac 1{|u|}u^T,$
\begin{equation}
\int_\Sigma e^{-\frac{X^2}4}\frac 1{|u|}\left(2k -|u|^2\right) dV\ge 0.
\end{equation}
But the assumption that $\Sigma$ is in $\overline{E^{k+1}(\sqrt{2k})}\times \mathbb R^{m-k-1}$ means
$$|u|^2-2k\ge 0$$
 Therefore, $|u|^2-2k=0,$
 i.e. $\Sigma\subset S^{k}(\sqrt{2k})\times \mathbb R^{m-k-1}.$

\begin{rem} 
If $n=m-1,$ then $\Sigma=S^{k}(\sqrt{2k})\times \mathbb R^{n-k}$ (\cite{caes}, Theorem 1.3).
\end{rem}
With the same arguments as in the proof of Theorem \ref{theovizh1}, we have the following theorem (see Colollary 1, \cite{vizh} for the case of codimension 1).
\begin{theorem} \label{theovizh2}
\begin{enumerate}
\item If the self-shrinker $\Sigma^n$ lies inside the closed cylinder 
$$\overline{B^{k+1}(a,\sqrt{2(n-p)+|a|^2}})\times\mathbb R^p,$$
where $a\in \mathbb R^{k+1},$ then $\Sigma\subset S^{k}(a,\sqrt{2(n-p)+|a|^2})\times\mathbb R^p.$ Moreover, if $n=m-1,$ then $\Sigma=S^{k}(\sqrt{2k})\times\mathbb R^{n-k}.$
\item The self-shrinker cannot lie outside the closed cylinder
$$\overline{B^{k+1}(a,\sqrt{2(k+1)+|a|^2}})\times\mathbb R^p,$$
for any vector $a$ in $\mathbb R^{k+1}.$
\end{enumerate}
\end{theorem}


\end{document}